\documentclass[12pt,fullpage]{article}

\title{\bf Van der Waerden/Schrijver-Valiant like Conjectures and Stable (aka Hyperbolic) Homogeneous Polynomials:\\
One Theorem for all}
\author{Leonid Gurvits\\ 
\small Los Alamos National Laboratory\\[-0.8ex]
\small \texttt{gurvits@lanl.gov}}

\begin{document}
\maketitle

\begin{abstract}
Let $p$ be a homogeneous polynomial of degree $n$ in $n$ variables, $p(z_1,\dots,z_n) = p(Z)$, $Z \in C^{n}$. 
We call such a polynomial $p$ {\bf H-Stable} if $p(z_1,\dots,z_n) \neq 0$ provided the real parts
$Re(z_i) > 0, 1 \leq i \leq n$. This notion from {\it Control Theory} is closely related
to the notion of {\it Hyperbolicity} used intensively in the {\it PDE} theory.

The main theorem in this paper states that if 
$p(x_1,\dots,x_n)$ is a homogeneous {\bf H-Stable} polynomial of degree $n$ with
nonnegative coefficients; 
$deg_{p}(i)$ is the maximum degree of the variable $x_i$, $C_i = \min(deg_{p}(i),i)$ and
$$
Cap(p) = \inf_{x_i > 0, 1 \leq i \leq n} \frac{p(x_1,\dots,x_n)}{x_1 \cdots x_n}  
$$
then the following inequality holds
$$
\frac{\partial^n}{\partial x_1\dots \partial x_n} p(0,\dots,0) \geq Cap(p) \prod_{2 \leq i \leq n} \left(\frac{C_i -1}{C_i} \right)^{C_{i}-1}. 
$$

This inequality is a vast (and unifying) generalization of the Van der Waerden conjecture on the permanents of doubly stochastic matrices as well as
the Schrijver-Valiant conjecture on the number of perfect matchings in $k$-regular bipartite graphs.
These two famous results correspond to the {\bf H-Stable} polynomials which are products of linear forms.   

Our proof is relatively simple and ``noncomputational''; it 
uses just very basic properties of complex numbers and the AM/GM inequality.
\end{abstract}

%  THEOREM-LIKE ENVIRONMENTS
 
\newtheorem{THEOREM}{Theorem}[section]
\newenvironment{theorem}{\begin{THEOREM} \hspace{-.85em} {\bf:} 
}%
                        {\end{THEOREM}}
\newtheorem{LEMMA}[THEOREM]{Lemma}
\newenvironment{lemma}{\begin{LEMMA} \hspace{-.85em} {\bf:} }%
                      {\end{LEMMA}}
\newtheorem{COROLLARY}[THEOREM]{Corollary}
\newenvironment{corollary}{\begin{COROLLARY} \hspace{-.85em} {\bf 
:} }%
                          {\end{COROLLARY}}
\newtheorem{PROPOSITION}[THEOREM]{Proposition}
\newenvironment{proposition}{\begin{PROPOSITION} \hspace{-.85em} 
{\bf:} }%
                            {\end{PROPOSITION}}
\newtheorem{DEFINITION}[THEOREM]{Definition}
\newenvironment{definition}{\begin{DEFINITION} \hspace{-.85em} {\bf 
:} \rm}%
                            {\end{DEFINITION}}
\newtheorem{EXAMPLE}[THEOREM]{Example}
\newenvironment{example}{\begin{EXAMPLE} \hspace{-.85em} {\bf:} 
\rm}%
                            {\end{EXAMPLE}}
\newtheorem{CONJECTURE}[THEOREM]{Conjecture}
\newenvironment{conjecture}{\begin{CONJECTURE} \hspace{-.85em} 
{\bf:} \rm}%
                            {\end{CONJECTURE}}
\newtheorem{PROBLEM}[THEOREM]{Problem}
\newenvironment{problem}{\begin{PROBLEM} \hspace{-.85em} {\bf:} 
\rm}%
                            {\end{PROBLEM}}
\newtheorem{QUESTION}[THEOREM]{Question}
\newenvironment{question}{\begin{QUESTION} \hspace{-.85em} {\bf:} 
\rm}%
                            {\end{QUESTION}}
\newtheorem{REMARK}[THEOREM]{Remark}
\newenvironment{remark}{\begin{REMARK} \hspace{-.85em} {\bf:} 
\rm}%
                            {\end{REMARK}}
\newtheorem{FACT}[THEOREM]{Fact}
\newenvironment{fact}{\begin{FACT} \hspace{-.85em} {\bf:} 
\rm}%                     
		            {\end{FACT}}

%\newenvironment{proof}{\noindent {\bf Proof:} \hspace{.677em}}%
%                      {}
 
%theorem
\newcommand{\thm}{\begin{theorem}}
%lemma
\newcommand{\lem}{\begin{lemma}}
%proposition
\newcommand{\pro}{\begin{proposition}}
%definition
\newcommand{\dfn}{\begin{definition}}
%remark
\newcommand{\rem}{\begin{remark}}
%example
\newcommand{\xam}{\begin{example}}
%conjecture
\newcommand{\cnj}{\begin{conjecture}}
%problem
\newcommand{\prb}{\begin{problem}}
%question
\newcommand{\que}{\begin{question}}
%corollary
\newcommand{\cor}{\begin{corollary}}
%fact
\newcommand{\fac}{\begin{fact}}

%proof
\newcommand{\prf}{\noindent{\bf Proof:} }
%end theorem
\newcommand{\ethm}{\end{theorem}}
%end lemma
\newcommand{\elem}{\end{lemma}}
%end proposition
\newcommand{\epro}{\end{proposition}}
%end definition
\newcommand{\edfn}{\bbox\end{definition}}
%end remark
\newcommand{\erem}{\bbox\end{remark}}
%end example
\newcommand{\exam}{\bbox\end{example}}
%end conjecture
\newcommand{\ecnj}{\bbox\end{conjecture}}
%end problem
\newcommand{\eprb}{\bbox\end{problem}}
%end question
\newcommand{\eque}{\bbox\end{question}}
%end corollary
\newcommand{\ecor}{\end{corollary}}
%end fact
\newcommand{\efac}{\end{fact}}
%end proof
\newcommand{\eprf}{\bbox}
%begin equation
\newcommand{\beqn}{\begin{equation}}
%end equation
\newcommand{\eeqn}{\end{equation}}
% white box
\newcommand{\wbox}{\mbox{$\sqcap$\llap{$\sqcup$}}}
%black box
\newcommand{\bbox}{\vrule height7pt width4pt depth1pt}
%\newcommand{\qed}{\bbox}

%right arrow
\newcommand{\rarrow}{\rightarrow}
%left arrow
\newcommand{\larrow}{\leftarrow}
%gradient
\newcommand{\grad}{\bigtriangledown}
\overfullrule=0pt
\def\setof#1{\lbrace #1 \rbrace}

\section{The permanent, the mixed discriminant, the Van Der Waerden conjecture(s) and homogeneous polynomials}
Recall that an $n \times n$ matrix $A$ is called doubly stochastic if it is nonnegative entry-wise
and its every column and row sum to one. The set of $n \times n$ doubly stochastic
matrices is denoted by $\Omega_{n}$. Let $\Lambda(k,n)$ denote the set of
$n \times n$ matrices with nonnegative integer entries and row and column sums all
equal to $k$. We define the following subset of rational doubly stochastic matrices:
$\Omega_{k,n} = \{ k^{-1} A: A \in \Lambda(k,n) \}$.
In a 1989 paper \cite{bapat}  R.B. Bapat defined the set 
$D_{n}$ of doubly stochastic $n$-tuples of $n \times n$ matrices.\\  
An $n$-tuple ${ \bf A }=(A_{1},
\dots, A_{n})$ belongs to $D_{n}$ iff 
$A_{i} \succeq 0$, i.e. $A_{i}$ is a positive
semi-definite matrix, $1 \leq i \leq n$;
$tr A_{i} = 1$ for $1 \leq i \leq n$;
$\sum^{n}_{i=1} A_{i} = I$, where $I$, as usual,
stands for the identity matrix.  
Recall that the permanent of a square matrix A is defined by
$$ per(A) = \sum_{\sigma \in S_{n}} \prod^{n}_{i=1} A(i,
\sigma(i)). $$
Let us consider an $n$-tuple
${ \bf A } = (A_{1}, A_{2},\dots  A_{n})$, where $A_{i} = (A_{i}(k,l):
1 \leq k, l \leq n)$ is a complex $n \times n$ matrix $(1 \leq i
\leq n)$.  Then 
$$
Det_{{\bf A}}(t_1,\dots,t_n) = \det (\sum_{1 \leq i \leq n} t_{i}A_{i})
$$ 
is a homogeneous
polynomial of degree $n$ in $t_{1},t_{2}, \dots, t_{n}$.  The number 
\begin{equation} \label{MD}
D({ \bf A }):=D(A_{1}, A_{2}, \dots, A_{n}) =
\frac{\partial^{n}}{\partial t_{1} \cdots \partial t_{n}} Det_{{\bf A}}(0,\dots,0)
\end{equation}
is called the mixed discriminant of $A_{1}, A_{2}, \dots,
A_{n}$.\\[1ex]
The mixed discriminant is just another name, introduced by A.D. Alexandrov, for $3$-dimensional Pascal's hyperdeterminant. 
The permanent is a particular (diagonal) case of the mixed discriminant. I.e.
define the following homogeneous polynomial 
\beqn \label{Prod}
Prod_{A}(t_1,\dots,t_n) = \prod_{1 \leq i \leq n} \sum_{1 \leq j \leq n} A(i,j)t_j.
\eeqn 
Then the next identity holds:
\beqn \label{PR}
per(A) = \frac{\partial^{n}}{\partial t_{1},\dots, \partial t_{n}} Prod_{A}(0,\dots,0).
\eeqn

Let us recall two famous results and one recent result by the author.

\bigbreak

\begin{enumerate}
\item {\bf Van der Waerden Conjecture}\\
The famous Van der Waerden Conjecture \cite{minc} states that
$$
\min_{A \in \Omega_{n}} per(A) = \frac{n!}{n^{n}} =: vdw(n)
\quad\hbox{ \bf (VDW-bound)}
$$
and the minimum is
attained uniquely at the matrix $J_{n}$ in which every entry equals
$\frac{1}{n}$. The Van der Waerden Conjecture was posed in 1926
and proved in 1981: D.I. Falikman proved
in \cite{fal} the lower bound
$\frac{n!}{n^{n}}$; 
the full conjecture, i.e. the uniqueness part, was proved by G.P. Egorychev in \cite{ego}.

\item {\bf Schrijver-Valiant Conjecture}\\
Define 
%\begin{align*}
%\lambda(k,n) &= \min \{per(A): A \in \Omega_{k,n}\} 
%= k^{-n} \min \{per(A): A \in \Lambda(k,n)\}; \\
%\theta(k) &= \lim_{n \rightarrow \infty}(\lambda(k,n))^{\frac{1}{n}}.
%\end{align*}
$$
\lambda(k,n)= \min \{per(A): A \in \Omega_{k,n}\} = k^{-n} \min \{per(A): A \in \Lambda(k,n)\};
$$
$$
\theta(k) = \lim_{n \rightarrow \infty}(\lambda(k,n))^{\frac{1}{n}}.
$$

It was proved in \cite{schr-val} that, using our notations, $\theta(k) \leq G(k) =: (\frac{k-1}{k})^{k-1}$ and conjectured
that $\theta(k) = G(k)$. 
Though the case of $k=3$ was proved by M. Voorhoeve in 1979 \cite{vor}, this conjecture was settled
only in 1998 \cite{schr} (17 years after the published proof of the Van der Waerden Conjecture). The main result of \cite{schr}
is the remarkable {\bf (Schrijver-bound)}:\\
\beqn \label{S-B}
\min \{per(A): A \in \Omega_{k,n}\} \geq \left(\frac{k-1}{k} \right)^{(k-1)n} 
\eeqn
The proof of {\bf (Schrijver-bound)} in \cite{schr} is, in the words of its author, ``highly complicated''.
\rem
The dynamics of research which led to {\bf (Schrijver-bound)}  is quite fascinating. If $k=2$ then $min_{A \in \Lambda(2,n)} per(A) = 2$.
Erdos and Renyi conjectured in 1968 paper that $3$-regular case already has exponential growth: 
$$
\min_{A \in \Lambda(3,n)} per(A) \geq a^n, a > 1.
$$
This conjecture is implied by {\bf (VDW-bound)}, this connection was another important motivation
for the Van der Waerden Conjecture. The Erdos-Renyi conjecture was answered by M. Voorhoeve in 1979 \cite{vor}:
\beqn \label{Voor} 
\min_{A \in \Lambda(3,n)} per(A) \geq 6 \left(\frac{4}{3} \right)^{n-3}.
\eeqn
Amazingly, the Voorhoeve's bound (\ref{Voor}) is asymptotically sharp and the proof of this fact is probabilistic.
In 1981 paper \cite{schr-val}, A.Schrijver and W.G.Valiant found a sequence $\mu_{k,n}$ of probabilistic distributions
on $\Lambda(k,n)$ such that 
\beqn \label{81}
\lim_{n \rightarrow \infty} \left(\min_{A \smash{\in \Lambda(k,n)}} per(A)\right)^{\frac{1}{n}} \leq    \lim_{n \rightarrow \infty} \left(E_{\mu_{k,n}} per(A) \right)^{\frac{1}{n}} = k \left(\frac{k-1}{k} \right)^{k-1}
\eeqn
({\it I.M. Wanless recently extended in \cite{wan} the upper bound (\ref{81}) to the boolean matrices in $\Lambda(k,n)$.)}\\

It follows from the Voorhoeve's bound (\ref{Voor}) that 
$$
\lim_{n \rightarrow \infty} \left(E_{\mu_{k,n}} per(A)\right)^{\frac{1}{n}} = \lim_{n \rightarrow \infty} \left(\min_{A\smash{\in \Lambda(k,n)}} per(A) \right)^{\frac{1}{n}} \quad \mbox{for}
\quad k=2,3.
$$
This was the rather bald intuition that gave rise to the {\it Schrijver-Valiant} 1981 conjecture.\\
The number $k \left(\frac{k-1}{k} \right)^{k-1}$ in {\it Schrijver-Valiant} conjecture came up via
combinatorics followed by the standard Stirling's formula manipulations. On the other
hand $G(k) = (\frac{k-1}{k})^{k-1} = \frac{vdw(k)}{vdw(k-1)}$.\\
\erem

\item {\bf Bapat's Conjecture (Van der Waerden Conjecture for mixed discriminants) }\\
One of the problems posed in \cite{bapat} is to determine the minimum
of mixed discriminants of doubly stochastic tuples:
$ min_{A \in D_{n}} D(A) = ?$ \\ 
Quite naturally, R.V.Bapat conjectured that
$ \min_{A \in D_{n}} D(A) = \frac{n!}{n^{n}}$ {\bf (Bapat-bound)}
and that it is attained uniquely at ${\bf J}_{n} =:(\frac{1}{n} I,\dots,
\frac{1}{n} I)$.\\
In \cite{bapat} this conjecture was formulated for real matrices. The author
proved it \cite{gur} even for the complex case, i.e. when matrices
$A_{i}$ above are complex positive semidefinite and, thus,
hermitian.
\end{enumerate}

\subsection{The Ultimate Unification (and Simplification)}
Falikman/Egorychev proofs of the Van Der Waerden conjecture as well our proof of Bapat's conjecture are based on the Alexandrov inequalities
for mixed discriminants \cite{Al38} and some optimization theory, which is rather advanced in the case of the Bapat's conjecture. 
They all rely heavily on the matrix structure and essentially of non-inductive nature.\\ 
{\it(D. I. Falikman independently rediscovered in \cite{fal} the diagonal case of the Alexandrov inequalities and used a
clever penalty functional. The very short paper \cite{fal} is supremely original, it cites only three references and uses
none of them.)}\\
The Schrijver's proof
has nothing in common with these analytic proofs; it is based on the finely tuned combinatorial arguments and
multi-level induction. It heavily relies on the fact that the entries of matrices $A \in \Lambda(k,n)$ are integers.

\medskip
 
The main result of this paper is one, easily stated and proved by easy induction, theorem which unifies, generalizes and, in the case
of {\bf (Schrijver-bound)}, improves the results described above. This theorem is formulated in terms of the mixed
derivative $\frac{\partial^n}{\partial x_1\dots \partial x_n} p(0,\dots,0)$ (rewind to the formula (\ref{PR}))
of {\bf H-Stable}  (or positive hyperbolic) homogeneous polynomials $p$.\\
The next two completely self-contained sections introduce the basics of stable homogeneous polynomials and proofs of the theorem and its corollaries.
We have tried to simplify everything to the undergraduate level, making the paper longer than a dry technical note of 4-5 pages.
Our proof of the uniqueness in the generalized Van der Waerden Conjecture is a bit more involved,
as it uses Garding's result on the convexity of the hyperbolic cone.

\section{Homogeneous Polynomials}
The next definition introduces key notations and notions.
\dfn
\begin{enumerate}
\item
The linear space of homogeneous polynomials with real (complex) coefficients
of degree $n$ and in $m$ variables is denoted $Hom_{R}(m,n)$ ($Hom_{C}(m,n)$).\\
We denote as $Hom_{+}(m,n)$ ($Hom_{++}(n,m)$) the closed convex cone
of polynomials $p \in Hom_{R}(m,n)$ with nonnegative (positive) coefficients.
\item
For a polynomial $p \in Hom_{+}(n,n)$ we define its {\bf Capacity} as
\beqn \label{cap}
Cap(p) = \inf_{x_i > 0, \prod_{1 \leq i \leq n }x_i = 1} p(x_1,\dots,x_n) = 
\inf_{x_i > 0} \frac{p(x_1,\dots,x_n)}{ \prod_{1 \leq i \leq n }x_i }.
\eeqn
\item
Consider a polynomial $p \in Hom_{C}(m,n)$, 
$$
p(x_1,\dots,x_m) = \sum_{(r_1,\dots,r_m)} a_{r_1,\dots,r_m} \prod_{1 \leq i \leq m } x_i^{r_{i}}. 
$$

We define
$Rank_{p}(S)$ as the maximal joint degree  attained on the subset\\ 
$S \subset \{1,\dots,m\}$:
\beqn \label{Rank-def}
Rank_{p}(S) = \max_{a_{r_1,\dots,r_m} \neq 0} \sum_{j \in S} r_j.
\eeqn

If $S = \{i\}$ is a singleton, we define $deg_{p}(i) = Rank_{p}(S)$.
\item
Let $p \in Hom_{+}(n,n)$, 
$$
p(x_1,\dots,x_n) = \sum_{r_1+ \cdots +r_n = n} a_{r_1,\dots,r_n} \prod_{1 \leq i \leq n } x_i^{r_{i}}. 
$$
Such a homogeneous polynomial $p$ with nonnegative coefficients is called doubly-stochastic if
$$
\frac{\partial}{\partial x_i} p(1,1,\dots,1) = 1: 1 \leq i \leq n.
$$
In other words, $p \in Hom_{+}(n,n)$ is doubly-stochastic if
\beqn \label{DS}
\sum_{r_1+ \cdots +r_n = n} a_{r_1,\dots,r_n} r_j = 1: 1 \leq j \leq n.
\eeqn
It follows from the Euler's identity that $p(1,1,\dots,1) = 1$:
\beqn \label{DS1}
\sum_{r_1+ \cdots +r_n = n} a_{r_1,\dots,r_n} = 1
\eeqn
Using the concavity of the logarithm on $R_{++}$ we get that
$$
\log \left(p(x_1,\dots,x_n) \right) \geq \sum_{r_1+ \cdots +r_n = n} a_{r_1,\dots,r_n} \sum_{1 \leq i \leq n} r_i \log(x_i) = \log(x_1 \cdots x_n).
$$
Therefore 
\fac \label{DS-Cap}
If $p \in Hom_{+}(n,n)$ is doubly-stochastic then $Cap(p) =1$.
\efac
\item
A polynomial $p \in Hom_{C}(m,n)$ is called {\bf H-Stable} if $p(Z) \neq 0$ provided
$Re(Z) > 0$; is called {\bf H-SStable} if $p(Z) \neq 0$ provided
$Re(Z) \geq 0$ and $\sum_{1 \leq i \leq m} Re(z_i) > 0$.\\
{\it We coined the term ``{\bf H-Stable}'' to stress two things: Homogeniety and Hurwitz' stability.
Other terms are used in the same context: {\bf Wide Sense Stable} in \cite{khar},
{\bf Half-Plane Property} in \cite{half}.}
\item
We define
\beqn \label{const}
vdw(i) = \frac{i!}{i^i}; G(i) = \frac{vdw(i)}{vdw(i-1)} = \left(\frac{i-1}{i} \right)^{i-1}, i > 1; G(1) = 1.
\eeqn
Notice that $vdw(i)$ as well  as $G(i)$ are strictly decreasing sequences.

\end{enumerate}
\edfn

\xam
\begin{enumerate}
\item
Let $p \in Hom_{+}(2,2), p(x_1,x_2) = \frac{A}{2} x_1^{2} + C x_1 x_2 + \frac{B}{2} x_2^{2}; A,B,C \geq 0$.
Then\\
$Cap(p) = C + \sqrt{AB}$ and the polynomial $p$ is {\bf H-Stable} iff $C \geq \sqrt{AB}$.
\item Let $A \in \Omega_{n}$ be a doubly stochastic matrix. Then the polynomial $Prod_{A}$ is doubly-stochastic.
Therefore $Cap(Prod_{A}) =1$. In the same way, if ${ \bf A } \in D_n$ is a doubly stochastic $n$-tuple then
the polynomial $Det_{{ \bf A }}$ is doubly-stochastic and
$Cap(Det_{{ \bf A }}) = 1$.
\item
Let ${ \bf A } = (A_{1}, A_{2},\dots  A_{m})$ be an $m$-tuple of PSD hermitian $n \times n$ matrices, and
$\sum_{1 \leq i \leq m} A_i \succ 0$ (the sum is positive-definite). Then the determinantal polynomial
$Det_{{\bf A}}(t_1,\dots,t_m) = \det (\sum_{1 \leq i \leq m} t_{i}A_{i})$ is {\bf H-Stable} and
\beqn \label{rank}
Rank_{Det_{{\bf A}}}(S) = Rank(\sum_{i \in S} A_{i}).
\eeqn

\end{enumerate}
\exam

The main result in this paper is the following Theorem.

\thm \label{main}
Let $p \in Hom_{+}(n,n)$ be {\bf H-Stable} polynomial. Then the
following inequality holds
\beqn \label{SCHR}
\frac{\partial^n}{\partial x_1\dots \partial x_n} p(0,\dots,0) \geq \prod_{2 \leq i \leq n} G \big(\min(i,deg_{p}(i)) \big) Cap(p).
\eeqn
\ethm 

Note that 
$$
\prod_{2 \leq i \leq n} G \big(\min(i,deg_{p}(i)) \big) \geq \prod_{2 \leq i \leq n} G(i) = vdw(n),
$$
which gives the next generalized {\bf Van Der Waerden Inequality}:

\cor \label{Waer}
Let $p \in Hom_{+}(n,n)$ be {\bf H-Stable} polynomial. Then
\beqn \label{Waer-Ineq}
\frac{\partial^n}{\partial x_1\dots \partial x_n} p(0,\dots,0) \geq \frac{n!}{n^n} Cap(p).
\eeqn
\ecor
{\it Corollary (\ref{Waer}) was conjectured by the author in \cite{newhyp}, where it was proved
that\\
 $\frac{\partial^n}{\partial x_1\dots \partial x_n} p(0,\dots,0) \geq C(n) Cap(p)$ for some constant $C(n)$.}\\

\subsection{Three Conjectures/Inequalities}
The fundamental nature of Theorem (\ref{main}) is illustrated in the following Example.
\xam
\begin{enumerate}
\item
Let $A \in \Omega_{n}$ be $n \times n$ doubly stochastic matrix.
It is easy to show that the polynomial $Prod_{A}$ is {\bf H-Stable} and doubly-stochastic.
Therefore $Cap(Prod_{A}) = 1$. Applying Corollary (\ref{Waer}) we get
the celebrated Falikman's result \cite{fal}:
$$
\min_{A \in \Omega_{n}} per(A) = \frac{n!}{n^n}.
$$
({\it The complementary uniqueness statement for Corollary (\ref{Waer}) will be considered in Section(\ref{kon}).})\\
\item
Let $(A_{1}, \dots, A_{n}) ={ \bf A } \in D_{n}$ be a doubly stochastic $n$-tuple.
Then the determinantal polynomial $Det_{{\bf A }}$ is {\bf H-Stable} and doubly-stochastic.
Thus $Cap(Det_{{\bf A }}) = 1$ and we get the {\bf (Bapat-bound)}, proved by the author:
$$
\min_{ { \bf A } \in D_{n}} D({ \bf A }) = \frac{n!}{n^n}.
$$ 

\item
Important for what follows is the next observation, which is a diagonal case of (\ref{rank}): \\
{\bf $deg_{Prod_{A}}(j)$ is equal to the number
of nonzero entries in the $j$th column of the matrix $A$.}\\ 
The next Corrolary combines this observation with Theorem(\ref{main}).
\cor \label{Improv}
\begin{enumerate} 
\item
Let $C_j$ be the number
of nonzero entries in the $j$th column of $A$, where $A$ is an $n \times n$ matrix with non-negative {\bf real} entries.
Then
\beqn \label{per-gen}
per(A) \geq \prod_{2 \leq j \leq n} G \left(\min(j,C_{j}) \right) Cap(Prod_{A}).
\eeqn
\item
Suppose that $C_j \leq k: k+1 \leq j \leq n$. Then
\beqn \label{schr-gen}
per(A) \geq \left(\Big(\frac{k-1}{k}\Big)^{k-1} \right)^{n-k} \frac{k!}{k^k} Cap(Prod_{A}).
\eeqn
\end{enumerate}
\ecor

Let $\Lambda(k,n)$ denote the set of
$n \times n$ matrices with nonnegative integer entries and row and column sums all
equal to $k$. The matrices in $\Lambda(k,n)$ correspond to the $k$-regular bipartite graphs
with multiple edges.\\ 
Recall the {\bf (Schrijver-bound)}:
\beqn \label{star}
\min_{A \in \Lambda(k,n)} per(A) \geq k^n G(k)^{n} = \left(\frac{(k-1)^{k-1}}{k^{k-2}} \right)^n.
\eeqn
The Falikman's inequality gives that 
$$
\min_{A \in \Lambda(k,n)} per(A) \geq k^n vdw(n) > k^n G(k)^n \quad \mbox{if} \quad k \geq n.
$$ 
Therefore the inequality (\ref{star}) is interesting only if $k < n$.\\
Note that if $A \in \Lambda(k,n), k < n$ then all columns of $A$ have at most $k$
nonzero entries.\\
If $A \in \Lambda(k,n)$ then the matrix $ \frac{1}{k} A \in \Omega_{n}$,  thus $Cap(Prod_{A}) = k^n$.
As we observed above, $deg_{Prod_{A}}(j) \leq k$. Applying the inequality (\ref{schr-gen}) to the polynomial
$Prod_{A}$ we get for $k < n$ an improved {\bf (Schrijver-bound)}: 
\beqn \label{new-schr}
\min_{A \in \Lambda(k,n)} per(A) 
\geq k^n \left(\Big(\frac{k-1}{k}\Big)^{k-1} \right)^{n-k} \frac{k!}{k^k} 
> \left(\frac{(k-1)^{k-1}}{k^{k-2}} \right)^n.
\eeqn
Interestingly, the inequality (\ref{new-schr}) recovers for $k=3$ the Voorhoeve's inequality (\ref{Voor}). 
\item
The inequality (\ref{per-gen}) is sharp if $C_i =\dots =C_{n-1}= n; C_{n} = k: 1 < k \leq n-1$.
To see this, consider the doubly stochastic matrix 

\beqn
D = \left( \begin{array}{cccc}
		  a &\dots & a & b \\
		  . &\dots  & . & . \\
		  a &\dots & a & b \\
		  c &\dots & c & 0 \\
		  . &\dots  & . & . \\
		  c &\dots & c & 0 \end{array} \right); a = \frac{1-b}{n-1} = \frac{k-1}{k(n-1)}, b =\frac{1}{k}, c = \frac{1}{n-1},
\eeqn

and the associated polynomial 
$$
Prod_{D}(x_1,\dots,x_n) = \left((\sum_{1 \leq i \leq n-1} ax_i) + bx_n \right)^{k} (\sum_{1 \leq i \leq n-1} cx_i)^{n-k}.
$$
Since the matrix $D$ is doubly stochastic, $Cap(Prod_{D})=1$. Direct inspection shows that 
$$
per(D) = (n-1)! (kb) a^{k-1} c^{n-k} = G(k) \frac{(n-1)!}{(n-1)^{n-1}}.
$$
Which gives the equality
$$
per(D) = Cap(Prod_{D}) \prod_{2 \leq j \leq n} G \left(\min(j,C_{j}) \right).
$$
It follows that $\min\{ per(A): A \in \Omega_{n}^{(0)} \} = \frac{(n-1)!}{(n-1)^{n-1}} \left(\frac{n-2}{n-1} \right)^{n-2}$,
where $\Omega_{n}^{(0)}$ is the set of $n \times n$ doubly stochastic matrices with at least one zero entry.
\end{enumerate} 
\exam
\subsection{The Main Idea}
Let $p \in Hom_{+}(n,n)$. Define the following polynomials $q_{i} \in Hom_{+}(i,i)$:
$$
q_{n} = p, q_{i}(x_1,\dots,x_i) = \frac{\partial^{n-i}}{\partial x_{i+1}\dots \partial x_n}p(x_1,\dots,x_i,0,\dots,0); 1 \leq i \leq n-1.
$$
\\
Notice that $q_{1}(x_{1}) = \frac{\partial^n}{\partial x_1\dots \partial x_n} p(0) x_1$ and
\beqn \label{deriv}
q_{2}(x_{1}, x_2) =\frac{\partial^n}{\partial x_1\dots \partial x_n} p(0) x_1x_2 + \frac{1}{2} \left(
\frac{\partial^n}{\partial x_1 \partial x_1\dots \partial x_n} p(0) x_1^{2} +  \frac{\partial^n}{\partial x_2 \partial x_2\dots \partial x_n} p(0) x_2^{2} \right).  
\eeqn
Therefore, $Cap(q_{1}) = \frac{\partial^n}{\partial x_1\dots \partial x_n} p(0)$ and
\beqn \label{deriv1}
Cap(q_{2}) = \frac{\partial^n}{\partial x_1\dots \partial x_n} p(0) +
\sqrt{\frac{\partial^n}{\partial x_1 \partial x_1\dots \partial x_n} p(0) \quad \frac{\partial^n}{\partial x_2 \partial x_2\dots \partial x_n}p(0)}.
\eeqn

Define the univariate polynomial $R(t)= p(x_1,\dots,x_{n-1},t)$. Then its derivative at zero is
\beqn \label{one-var}
R^{\prime}(0) =q_{n-1}(x_1,\dots,x_{n-1}).
\eeqn 

Another simple but important observation is the next inequality:
\beqn \label{monot}
deg_{q_{i}}(i) \leq \min \left(i,deg_{p}(i) \right) \Longleftrightarrow G \left(deg_{q_{i}}(i) \right) \geq G \left(\min(i,deg_{p}(i)) \right): 1 \leq i \leq n.
\eeqn 
Recall that $vdw(i) = \frac{i!}{i^i}$. Suppose that the next inequalities hold
\beqn \label{ind}
Cap(q_{i-1}) \geq Cap(q_{i}) \frac{vdw(i)}{vdw(i-1)} = Cap(q_{i}) G(i): 2 \leq i \leq n.
\eeqn

Or better, the next stronger ones hold
\beqn \label{indS}
Cap(q_{i-1}) \geq Cap(q_{i}) G \left(deg_{q_{i}}(i) \right): 2 \leq i \leq n,
\eeqn
where
\beqn
G(m) = \frac{vdw(m)}{vdw(m-1)} = \left(\frac{m-1}{m} \right)^{m-1}.
\eeqn
The next result, proved by the straigthforward induction, summarizes the main idea of our approach.

\thm
\begin{enumerate}
\item
If the inequalities (\ref{ind}) hold then the next generalized Van Der Waerden inequality holds:
\beqn \label{VDW}
\frac{\partial^n}{\partial x_1\dots \partial x_n} p(0,\dots,0) = Cap(q_{1}) \geq vdw(n)Cap(p).
\eeqn
In the same way, the next inequality holds for $Cap(q_{2})$:
\beqn \label{VDW1}
\frac{\partial^n}{\partial x_1\dots \partial x_n} p(0) +
\sqrt{\frac{\partial^n}{\partial x_1 \partial x_1\dots \partial x_n} p(0) \frac{\partial^n}{\partial x_2 \partial x_2\dots \partial x_n}p(0)} \geq 2 vdw(n)Cap(p).
\eeqn

\item
If the inequalities (\ref{indS}) hold then the next generalized {\bf (Schrijver-bound)} holds:
\beqn \label{SCHR1}
\frac{\partial^n}{\partial x_1\dots \partial x_n} p(0,\dots,0) = Cap(q_{1}) \geq Cap(p) \prod_{2 \leq i \leq n} G \big(\min(i,deg_{p}(i)) \big).
\eeqn
\end{enumerate}
\ethm

What is left is to prove that the inequalities (\ref{indS}) hold for {\bf H-Stable} polynomials.\\
We break the proof of this statement in two steps.
\begin{enumerate}
\item Prove that if $p \in Hom_{+}(n,n)$ is {\bf H-Stable} then $q_{n-1}$ is either zero or {\bf H-Stable}.
Using equation (\ref{one-var}), this implication follows from Gauss-Lukas Theorem. Gauss-Lukas Theorem
states that if $z_1,\dots,z_n \in C$ are the roots of an univariate polynomial $Q$ then
the roots of its derivative $Q^{\prime}$ belong to the convex hull $CO(\{z_1,\dots,z_n\})$.\\
{\it This step is, up to minor perturbation arguments, known. See, for instance, \cite{khov}. The result in
\cite{khov} is stated in terms of hyperbolic polynomials, see Remark (\ref{khov}) for the connection between
{\bf H-Stable} and hyperbolic polynomials. Our treatment, described in Section(\ref{shp}), is self-contained,
short and elementary.}
\item Prove that $Cap(q_{n-1}) \geq G(deg_{p}(n)) Cap(p)$. This inequality boils down
to the next inequality for the univariate polynomial $R$ from (\ref{one-var}):
$$
R^{\prime}(0) \geq G(deg(R)) \left(\inf_{t > 0} \frac{R(t)}{t} \right).
$$
We prove it using AM/GM inequality and the fact that the roots of the polynomial $R$ are real.
\end{enumerate}

{\it It is instructive to see what is going on in the ``permanental case'': we start with the
polynomial $Prod_{A}$ which is a product of nonnegative linear forms. The very first polynomial
in the induction, $q_{n-1}$, is not of this type in the generic case. I.e. there is no
one matrix/graph associated with $q_{n-1}$. We gave up the matrix structure but had won the game.} 

In the rest of the paper {\bf Facts} are statements which are quite simple and (most likely) known.
We included them having in mind the undergraduate student reader.
\section{Univariate Polynomials}
\pro \label{GL}
\begin{enumerate}
\item (Gauss-Lukas Theorem)\\
Let $R(z) = \sum_{0 \leq i \leq n} a_i z^i$ be a Hurwitz polynomial with complex coefficients,
i.e. all the roots of $R$ have negative real parts.\\
Then its derivative $R^{\prime}$ is Hurwitz.  
\item
Let $R(z) = \sum_{0 \leq i \leq n} a_i z^i$ be a Hurwitz polynomial
with real coefficients and $a_n > 0$. Then all the coefficients
are positive real numbers.
\end{enumerate}
\epro

\prf
\begin{enumerate}
\item
Recall that 
$$
\frac{R^{\prime}(z)}{R(z)} = \sum_{1 \leq j \leq n} \frac{1}{z - z_j}.
$$
Let $\mu$ be a root of $R^{\prime}$. Consider two cases.
First: $\mu$ is a root of $R$. Then clearly $Re(\mu) < 0$.
Second: $\mu$ is not a root of $R$. Then
$$
L =: \sum_{1 \leq j \leq n} \frac{1}{\mu - z_j} = 0.
$$
Suppose that $Re(\mu) \geq 0$. As $(a + ib)^{-1} = \frac{a - ib}{a^2 + b^2}$
we get that 
$$
Re \left(\frac{1}{\mu - z_j} \right) = \frac{Re (\mu) - Re(z_j)}{\left(Re(\mu) - Re(z_j) \right)^2 + \left(Im(\mu) - Im(z_j) \right)^2} > 0.
$$
Therefore $Re(L) > 0$ which leads to a contradiction. Thus $Re(\mu) < 0$ and the derivative $R^{\prime}$ is Hurwitz.
\item This part is easy and well known.
\end{enumerate}
\eprf

The next simple result binds together all the small pieces of our approach.

\lem \label{derest}
 Let $Q(t) = \sum_{0 \leq i \leq k} a_i t^i; a_k > 0, k \geq 2$ be a polynomial with non-negative coefficients
 and real (non-positive) roots. Define $C = \inf_{t>0} \frac{Q(t)}{t}$. Then the
 next inequlity holds:
 \beqn \label{hui}
 a_1 = Q^{\prime}(0) \geq \left(\frac{k-1}{k} \right)^{k-1} C.
 \eeqn
 The equality holds if and only if all the roots of $Q$ are equal negative numbers, i.e.
 $Q(t) = b (t+a)^{k}$ for some $a,b > 0$.
 \elem

 \prf 
 If $Q(0)=0$ then $Q^{\prime}(0) \geq C > \left(\frac{k-1}{k} \right)^{k-1} C $.\\
 Let $Q(0)> 0$. We then can assume WLOG that $Q(0) = 1$. In this case
 all the roots of $Q$ are negative real numbers. Thus
 $$
  Q(t):=\prod_{i=1}^k (a_{i}t + 1): a_i > 0, 1 \leq i \leq k,
 $$
 and $Q^{\prime}(0) = a_1+\dots +a_k$.\\
 
 Using the AM/GM inequality we get that
 \beqn \label{conc}
 Ct \leq Q(t) \leq P(t) =: \left(1 + \frac{Q'(0)}{k} t \right)^k, t \geq 0.
 \eeqn
 It follows from basic calculus that
 $$
 \inf_{t>0} \frac{P(t)}{t} = P(s) = Q'(0)\left(\frac{k}{k-1} \right)^{k-1}, \quad \mbox{where} \quad s =\frac{k}{Q'(0)(k-1)}.
 $$
 
 Therefore
 $$
 C \leq \inf_{t>0} \frac{P(t)}{t} = Q'(0)\left(\frac{k}{k-1} \right)^{k-1},
 $$
 which finally yields the desired inequality
 $$
 Q'(0) \geq \left(\frac{k-1}{k} \right)^{k-1} C, k \geq 2.
 $$
 
 It follows from the uniqueness condition in the AM/GM inequality that
 the equality in (\ref{hui}) holds if and only if $0 < a_1 =\dots = a_k $.
 \eprf 
 \rem The condition that the roots of $Q$ are real can be relaxed in several ways.
 For instance the statement of Lemma (\ref{derest}) holds for any map $f:R_{+} \rightarrow R_{+}$
 such that the derivative $f^{\prime}(0)$ exists and $f^{\frac{1}{k}}$ is concave.\\
 If such map is log-concave, i.e $\log (f)$ is concave, then $f^{\prime}(0) \geq \frac{1}{e} \inf_{t>0} \frac{f(t)}{t}$.\\
 Notice that the right inequality in (\ref{conc}) is essentially equivalent to the concavity\\
 of the function $\left(Q(t) \right)^{\frac{1}{k}}$ on $R_{+}$.\\
 It was shown in \cite{my-vdw} that the inequality (\ref{hui}) is equivalent to the {\bf (VDW-bound)} for
 doubly stochastic matrices $A \in \Omega_{n}: A = [a|b|\dots |b]$ with two distinct columns. 
 \erem

\section{Stable homogeneous polynomials} \label{shp}

\subsection{Basics}
\dfn
A polynomial $p \in Hom_{C}(m,n)$ is called {\bf H-Stable} if $p(Z) \neq 0$ provided
$Re(Z) > 0$; is called {\bf H-SStable} if $p(Z) \neq 0$ provided
$Re(Z) \geq 0$ and $\sum_{1 \leq i \leq m} Re(z_i) > 0$.
\edfn

\fac \label{approx}
Let $p \in Hom_{C}(m,n)$ be {\bf H-Stable} and $A$ is $m \times m$ matrix with nonnegative
real entries without zero rows. Then the polynomial
$p_{A}$, defined as $p_{A}(Z) = p(AZ)$ is also {\bf H-Stable}. If
all entries of $A$ are positive then $p_{A}$ is {\bf H-SStable}. 
\efac

\fac \label{rts}
Let $p \in Hom_{C}(m,n)$, $Y \in C^m, p(Y) \neq 0$. Define the following univariate
polynomial of degree $n$: 
$$
L_{X,Y}(t) = p(t Y - X) = p(Y) \prod_{1 \leq i \leq n} (t - \lambda_{i;Y}(X)): X \in C^m.
$$
Then 
\beqn \label{det}
\lambda_{i;Y}(bX + aY) = b \lambda_{i;Y}(X) + a; p(X) = p(Y) \prod_{1 \leq i \leq n} \lambda_{i;Y}(X).
\eeqn
\efac

The following simple result substantially simplifies the proofs below. Proposition (\ref{st-hyp}) connects
the notion of {\bf H-Stability} with the notion of {\bf Hyperbolicity}, see more on this connection
in Subsection(5.1).

\pro \label{st-hyp}
A polynomial $p \in Hom_{C}(m,n)$ is {\bf H-Stable} if and only if
$p(X) \neq 0: X \in R^{m}_{++}$ and the roots of univariate
polynomials $P(tX-Y): X,Y \in R^{m}_{++}$ are real positive numbers.
\epro

\prf
\begin{enumerate}
\item
Suppose that $p(X) \neq 0: X \in R^{m}_{++}$ and the roots of univariate
polynomials $p(tX-Y): X,Y \in R^{m}_{++}$ are real positive numbers. 
It follows from identities (\ref{det}) (shift $L \rightarrow L + aX > 0$) that the roots of $P(tX-L): X \in R^{m}_{++}, L \in R^m $ are real numbers.
We want to prove that
this property implies that $p \in Hom_{C}(m,n)$ is {\bf H-Stable}. Let $Z = Re(Z) + i Im(Z) \in C^m: Im(Z) \in R^m, 0 < Re(Z) \in R^{m}_{++}$.
If $p(Z) = 0$ then also $p(-i Re(Z) + Im(Z)) = 0$, which contradicts the real rootedness of $p(tX -Y): X > 0, Y \in R^m$.  
\item
Suppose that $p \in Hom_{C}(m,n)$ is {\bf H-Stable}. Let $X,Y \in R^{m}_{++}$ and $p(zX-Y) = 0, z = a + bi$. We need to prove that
$b = 0$ and $a > 0$. If $b \neq 0$ then $p(aX -Y + bi X) = (bi)^n p(X - b^{-1} i(aX - y)) \neq 0$ as 
the real part $Re(X - b^{-1} i(aX - y)) = X > 0$. Therefore $b = 0$. If $a \leq 0$ then $-(aX -Y) \in R^{m}_{++}$.
Which implies that $p(aX -Y) = (-1)^n p(-(aX -Y)) \neq 0$. Thus $a > 0$.
\end{enumerate}
\eprf

We will use the following corollaries:
\cor
If $Re(Z) \in R^{m}_{+}$ and a polynomial $p$ is {\bf H-Stable}  then
\beqn \label{max-prin}
|p(Z)| \geq |p \left(Re(Z) \right)|.
\eeqn
\ecor
\prf
Since $p$ is continuous on $C^m$ hence it is sufficient to assume that $Re(Z) \in R^{m}_{++}$.\\
It follows from identities (\ref{det}) that
\beqn \label{roots}
p(Z) = p \left(Re(Z) + i Im(Z) \right) = p \left(Re(Z)\right) \prod_{1 \leq j \leq n} (1+ i \lambda_j),
\eeqn
\\
where $(\lambda_1,\dots,\lambda_n)$ are the roots of the univariate polynomial
$p(t Re(Z) - Im(Z))$. Because $Re(Z) \in R^{m}_{++}$, all these roots are real numbers.\\
Therefore $|p(Z)| = |p(Re(Z))| \prod_{1 \leq j \leq n} |1+ i \lambda_j | \geq  |p(Re(Z))|$.
\eprf

\cor \label{two-vect}
Let $p \in Hom_{C}(m,n)$ be {\bf H-Stable}; $X,Y \in R^m$ and $0 < X+Y \in R^m_{++}$.
Then all the roots of the univariate polynomial equation $p(tX + Y)= 0$ are real numbers.
\ecor 
\prf
Let $p(tX + Y)= 0$, then also $p((t-1)X + (X+Y))=0$. Since $X+Y > 0$ hence $t-1 \neq  0$.
As the polynomial $p$ is homogeneous therefore $p(X + (1-t)^{-1} (X+Y)) = 0$.
It follows that $(1-t)^{-1}$ is real, thus $t$ is also a real number.
\eprf

\fac \label{Pos}
Let $p \in Hom_{C}(m,n)$ be {\bf H-SStable} ({\bf H-Stable}). Then for all $X \in R^{m}_{++}$  the coefficients of the polynomial
$q = \frac{p}{p(X)}$ are positive (nonnegative) real numbers. 
\efac

\prf
We prove first the case of {\bf H-SStable} polynomials.\\
Since $q(X) = 1$ we get from (\ref{det}) that  $q(Y)$ is a positive real number for  all vectors $Y \in R^{m}_{++}$.
Therefore, by a standard interpolation argument, the coefficients of $q$ are real. We will prove by induction the following equivalent
statement: {\it if $q \in Hom_{R}(m,n)$ is {\bf H-SStable} and $q(Y) > 0$ for all $Y \in R^{m}_{++}$ 
then the coefficients of $q$ are all positive.}
Write $q(t; Z) = \sum_{0 \leq i \leq n} t^{i} q_i(Z)$, where $Z \in C^{m-1}$,
the polynomials $q_i \in Hom_{R}(m-1,n-i)$, $0 \leq i \leq n-1$ and $q_{n}(Z)$ is a real number.
Let us fix the complex vector $Z$ such that $Re(Z) \in R^{m-1}_{+}$ and $Re(Z) \neq 0$. Since $q$ is {\bf H-SStable} hence
all roots of the univariate polynomial $q(t; Z)$ have negative real parts. Therefore, using the first part of Proposition (\ref{GL}),
we get that polynomials $q_i: 0 \leq i \leq n$ are all {\bf H-SStable}. Since the degree of $q$ is $n$
hence $q_{n}(Z)$ is a constant, $q_{n}(Z) = q(1; 0) > 0$. Using now the second part of Proposition (\ref{GL}),
we see that $q_{i}(Y) > 0$ for all $Y \in R^{m}_{++}$ and $0 \leq i \leq n$. Continuing this
process we will end up with either $m=1$ or $n=1$. Both those cases have positive coefficients.\\

Let $p \in Hom_{C}(m,n)$ be {\bf H-Stable} and $A > 0 $ is $m \times m$ matrix with positive entries such that
$AX =X$.
Then for all $\epsilon > 0$ the polynomials $q_{I + \epsilon A} \in Hom_{R}(m,n)$, defined as in Fact(\ref{approx}), are {\bf H-SStable} and
$\lim_{\epsilon \rightarrow 0}q_{I + \epsilon A} = q$. Therefore the coefficients of $q$ are nonnegative real numbers.  
\eprf

{}From now on we will deal only with the polynomials
with nonnegative coefficients.

\cor \label{limit}
Let $p_i \in Hom_{+}(m,n)$ be a sequence of {\bf H-Stable} polynomials and $p = \lim_{i \rightarrow \infty} p_i$.
Then $p$ is either zero or {\bf H-Stable}.
\ecor
{\it Some readers might recognize Corollary (\ref{limit}) as a particular case of A. Hurwitz's theorem on limits
of sequences of nowhere zero analytical functions. Our proof below is elementary.}\\ 

\prf
Suppose that $p$ is not zero. Since $p \in Hom_{+}(m,n)$ hence $p(x_1,\dots,x_m) > 0$ if $x_j > 0: 1 \leq j \leq m$.
As the polynomials $p_i$ are {\bf H-Stable} therefore $|p_{i}(Z)| \geq |p_{i} \left(Re(Z) \right)|: Re(Z) \in R_{++}^{m}$.
Taking the limits we get that $|p(Z)| \geq |p \left(Re(Z) \right)| > 0: Re(Z) \in R_{++}^{m}$, which means
that $p$ is {\bf H-Stable}. 
\eprf

\fac \label{Recur}
For a polynomial $p \in Hom_{C}(m,n)$ we define a polynomial\\
$q \in Hom_{C}(m-1,n-1)$ as
$$
q(x_1,\dots,x_{m-1}) = \frac{\partial}{\partial x_m}p(x_1,\dots,x_{m-1},0).
$$
Then the next two statements hold:
\begin{enumerate}
\item
Let $p \in Hom_{+}(m,n)$ be {\bf H-SStable}. Then
the polynomial $q$ is also {\bf H-SStable}.
\item
Let $p \in Hom_{+}(m,n)$ be {\bf H-Stable}. Then
the polynomial $q$ is either zero or {\bf H-Stable}.
\end{enumerate}
\efac

\prf
\begin{enumerate}
\item
Let $p \in Hom_{+}(m,n)$ be {\bf H-SStable} and consider
an univariate polynomial 
$$
R(z) = p(Y;z): z \in C, Y \in C^{m-1}.
$$
Suppose that $0 \neq Re(Y) \geq 0$. It follows from the definition
of {\bf H-SStability} that $R(z) \neq 0$ if $Re(z) \geq 0$.
In other words, the univariate polynomial $R$ is Hurwitz.
It follows from  Gauss-Lukas Theorem that
$$
q(Y) = R^{\prime}(0) \neq 0,
$$
which means that $q$ is {\bf H-SStable}.
\item
Let $p \in Hom_{+}(m,n)$ be {\bf H-Stable} and $q \neq 0$.
Take an $m \times m$ matrix $A > 0$.
Then the polynomial $p_{I + \epsilon A}$, 
$p_{I + \epsilon A}(Z) = p \left((I + \epsilon A)Z \right)$
is {\bf H-SStable} for all $\epsilon > 0$. Therefore, using the first part,
$q_{I + \epsilon A}$ is {\bf H-SStable}. Clearly
$\lim_{\epsilon \rightarrow 0} q_{I + \epsilon A} = q$.
Since $q \neq 0$, it follows from Corollary (\ref{limit}) that
$q$ is {\bf H-Stable}. 
\end{enumerate}
\eprf

\thm \label{vse}
Let $p \in Hom_{+}(n,n)$ be {\bf H-Stable},
and 
$$
q_{n-1}(x_1,\dots,x_{n-1}) = \frac{\partial}{\partial x_n}p(x_1,\dots,x_{n-1},0).
$$
Then
\beqn \label{Cap-In}
Cap(q_{n-1}) \geq Cap(p) G \left(deg_{p}(n) \right).
\eeqn
\ethm

\prf
We need to prove that
$$
\frac{\partial}{\partial x_n}p(x_1,\dots,x_{n-1},0) \geq Cap(p) G \left(deg_{p}(n) \right), x_1,\dots,x_{n-1} > 0, \prod_{1 \leq i \leq n-1} x_i =1.
$$ 
Fix a positive vector $(x_1,\dots,x_{n-1}),\prod_{1 \leq i \leq n-1} x_i =1$ and define, as in proof of Fact (\ref{Recur}),\\ 
the polynomial $R(t) = p(x_1,\dots,x_{n-1},t)$. It follows from Corollary(\ref{two-vect}) that all the roots
of $R$ are real. Since the coefficients of the polynomial $R$ are non-negative hence its roots are non-positive real numbers.
It follows from a definition of $Cap(p)$ that $R(t) \geq Cap(p) t$, therefore  
$$
\inf_{t > 0} \frac{R(t)}{t} \geq Cap(p).
$$
The degree of the polynomial $R$ is equal to $deg_{p}(n)$. It finally follows from Lemma(\ref{derest}) that
$$
q_{n-1}(x_1,\dots,x_{n-1},0) = R^{\prime}(0) \geq Cap(p)G \left(deg_{p}(n) \right).
$$
\eprf

\section{Uniqueness in Generalized Van Der Waerden Inequality} \label{kon}
\subsection{Hyperbolic Polynomials}
{\it The following concept  of hyperbolic polynomials arose from the theory of partial differential equations \cite{gar}, 
\cite{horm}. A recent paper \cite{ren} gives nice
and concise introduction to the area (with simplified proofs of the key theorems) and describes connections to
convex optimization.}

\bigbreak

\dfn
\begin{enumerate}
\item
A homogeneous polynomial $p: C^m \rightarrow C$ of degree $n$ ($p \in Hom_{C}(m,n)$) is called hyperbolic in the direction $e \in R^m$ 
(or $e$-hyperbolic) if $p(e) \neq 0$ and for each vector $X \in R^m$ the univariate (in $\lambda$) polynomial 
 $p(X - \lambda e)$   has exactly $n$ real roots counting their multiplicities. 
\item
Denote an ordered vector of roots of $p(x - \lambda e)$ as 
$$
\lambda_{e}(X) = (\lambda_{n}(X) \geq \lambda_{n-1}(X) \geq \dots  \lambda_{1}(X)).
$$ 
Call $X \in R^m$ $e$-positive ($e$-nonnegative) if
$\lambda_{1}(X) > 0$ ($\lambda_{n}(X) \geq 0$). 
We denote the closed set of $e$-nonnegative vectors as $N_{e}(p)$,
and the open set of $e$-positive vectors
as $C_{e}(p)$.
\end{enumerate}
\edfn

\rem \label{khov}
Proposition (\ref{st-hyp}) essentially says that a polynomial $p \in Hom_{C}(m,n)$ is {\bf H-Stable} iff
$p$ is hyperbolic in some direction $e \in R^m_{++}$ and the inclusion
$R^m_{+} \subset N_{e}(p)$ holds. If $p \in Hom_{C}(m,n)$ is {\bf H-SStable} then
any non-zero vector $0 \leq X \in R^{m}_{+}$ belongs to the (open) {\it hyperbolic cone} $C_{e}(p)$. 
\erem

We need the next fundamental fact due to L. Garding \cite{gar} (we recommend the very readable treatment in \cite{ren}):
\thm \label{gard}
Let $p \in Hom_{C}(m,n)$ be $e$-hyperbolic polynomial and $d \in C_{e}(p) \subset R^m$. Then $p$ is also $d$-hyperbolic and 
$C_{d}(p) = C_{e}(p), N_{d}(p) = N_{e}(p) $. Moreover cone $C_{e}(p)$, called {\it hyperbolic cone}, is convex. 
\ethm

\cor \label{lin}
\begin{enumerate}
\item
For any two vectors in the hyperbolic cone $d_1, d_2 \in C_{e}(p)$ the following set equality holds:
\beqn \label{same} 
N_{d_1}(p) \bigcap \left(-N_{d_1}(p) \right) = N_{d_2}(p) \bigcap \left(-N_{d_2}(p) \right) = Null_{p}.
\eeqn
Thus $Null_{p} \subset R^m$ is a linear subspace.
\item

\beqn \label{same-more}
Null_{p} = \{X \in R^m: p(Y + X) = p(Y) \quad \mbox{for all} \quad Y \in C^m \}
\eeqn
Let $Pr(Null_{p})$ be orthogonal projector on the linear subspace $Null_{p}$. 
It follows from (\ref{same-more}) that
\beqn \label{proj} 
p(Y) = p \left((I-Pr(Null_{p})) Y \right)
\eeqn
\item
Let $p \in Hom_{+}(n,n)$ be a doubly-stochastic {\bf H-Stable} polynomial. If $(y_1,\dots,y_n) \in Null_{p}$
then
\beqn \label{trace}
y_1+\dots +y_n = 0.
\eeqn

\end{enumerate}
\ecor

\bigbreak

\prf
\begin{enumerate}
\item It is well known and obvious that if $K$ is a convex cone in some linear space $L$ over reals then
the intersection $K \bigcap (-K)$ is a linear subspace of $L$.
\item
Let $T \in C_{e}(p)$ and $X \in Null_{p}$. Then all the roots of the equation $p(xT + X) =0$ are equal to zero.
Since $p(T) \neq 0$ and the polynomial $p \in Hom_{C}(m,n)$ is homogeneous, hence $p(xT + X) = x^n p(T)$. 
Therefore $p(T + X) = p(T)$ for all $T \in  C_{e}(p)$. As $C_{e}(p)$ is a non-empty open subset of $R^m$,
equality (\ref{same-more}) follows from the analyticity of $p$. 
\item Consider the vector of all ones $e = (1,\dots,1) \in R^n$ and a vector $Y = (y_1,\dots,y_n) \in Null_{p}$.
Then $d(t) = p(e + tY) = p(e)$ for all $t \in R$. Therefore
$$
0 = d^{\prime}(0) = \sum_{1 \leq i \leq n} y_i \frac{\partial}{\partial x_i} p(1,1,\dots,1)  = y_1+\dots +y_n.
$$

\end{enumerate}

\eprf

\xam \label{one-dim}
\begin{enumerate}
\item
Consider the power polynomial $q \in Hom_{+}(n,n)$, $q(x_1,\dots,x_n) = (a_1 x_1+\dots +a_n x_n)^n$.
If the non-zero vector ${\bf a} = (a_1,\dots,a_n) \in R^n_{+}$ then the power polynomial $q$ is {\bf H-Stable}.
The correspondind linear subspace 
$$
Null_{p} = {\bf a}^{\perp} =: \{(y_1,\dots,y_n) \in R^n: \sum_{1 \leq i \leq n} a_i y_i = 0 \}, \quad dim(Null_{p}) = n-1;
$$
and $Cap(p) = n^n \prod_{1 \leq i \leq n} a_i$. Therefore $Cap(q) \neq 0$ iff ${\bf a} \in R^n_{++}$.\\
It is easy to see that
$$
\frac{\partial^n}{\partial x_1\dots \partial x_n} (a_1 x_1+\dots +a_n x_n)^n = n! a_1\dots a_n.
$$
Therefore 
\beqn \label{extr}
\frac{\partial^n}{\partial x_1\dots \partial x_n} q(0,\dots,0) = Cap(q) \frac{n!}{n^n}.
\eeqn

If $dim(Null_{p}) = n-1$ and a polynomial $p \in Hom_{+}(n,n)$ is {\bf H-Stable} then $p(x_1,\dots,x_n) = (b_1 x_1+\dots +b_n x_n)^n$
for some non-zero vector ${\bf b} = (b_1,\dots,b_n) \in R^n_{+}$. The power polynomial $p$ is
doubly-stochastic iff $b_i = \frac{1}{n}, 1 \leq i \leq n$.
\item
Let $p \in Hom_{C}(m,n)$ be an $e$-hyperbolic polynomial, $D \in C_{e}(p) \subset R^m$ and $X \in R^m$.
Suppose that the univariate polynomial $R(t) = p(tD + X) = a (t + b)^n, b \in R$. Define the
next real vector $Y = -bD + X$. Then all the roots of the equation $p(Y - \lambda D) = 0$ are equal to zero.
Therefore $Y \in N_{D}(p) \bigcap \left(-N_{D}(p) \right) = Null_{p}$.

\end{enumerate}
\exam

\subsection{Uniqueness}
\dfn
We call a {\bf H-Stable} polynomial $p \in Hom_{+}(n,n)$ {\bf extremal} if\\
$Cap(p) > 0$ and\\
\beqn \label{EQ}
\frac{\partial^n}{\partial x_1\dots \partial x_n} p(0,\dots,0) = \frac{n!}{n^n} Cap(p).
\eeqn
\edfn

Our goal is the next theorem
\thm \label{uniq}
A {\bf H-Stable} polynomial $p \in Hom_{+}(n,n)$ is {\bf extremal} if and only if
$$
p(x_1,\dots,x_n) = (a_1 x_1 +\dots + a_n x_n)^{n}
$$
for some positive real numbers $a_1,\dots,a_n > 0$.\\
In other words, the equality (\ref{EQ}) holds iff $dim(Null_{p}) = n-1$.
\ethm 
{\it Notice that the ``if'' part is simple and follows from the equality(\ref{extr}).}

We collect the basic properties of {\bf extremal} polynomials in the next proposition.

\pro \label{ex-pr}
\begin{enumerate}
\item If a {\bf H-Stable} polynomial $p \in Hom_{+}(n,n)$ is {\bf extremal} then
all its coefficients are positive real numbers.
\item
Let ${\bf c} =(c_1,\dots,c_n) \in R_{++}^{n}$ and $p \in Hom_{+}(n,n)$. Define
the scaled polynomial $p_{{\bf c}}$ as $p_{{\bf c}}(x_1,\dots,x_n) = p(c_1 x_1,\dots,c_n x_n)$.
If $p \in Hom_{+}(n,n)$ is {\bf H-Stable} {\bf extremal} polynomial then also the scaled polynomial
$p_{{\bf c}}$ is.
\end{enumerate}
\epro

\prf
\begin{enumerate}
\item
Our goal in this step is to show that if equation (\ref{EQ}) holds then
the {\bf H-Stable} polynomial $p$ is, in fact, {\bf H-SStable} and therefore has all positive coefficients.

Since $G(2) \cdots G(n) = \frac{n!}{n^n}$ and the function $G$ is strictly decreasing on $[0, \infty)$, hence 
it follows from (\ref{SCHR}) that $deg_{p}(n) =n$. Since the inequality (\ref{SCHR}) is invariant
with respect to permutations of variables hence $deg_{p}(i) =n: 1 \leq i \leq n$. Which means that
$p(e_i) > 0: 1 \leq i \leq n$, where $\{e_1,\dots,e_n\}$ is the standard orthonormal basis in $R^n$.
Therefore the polynomial $p$ is {\bf H-SStable}. Thus its coefficients
are strictly positive real numbers and all non-zero vectors $Y \in R_{+}^n$ belong to its open
hyperbolic cone $C_{e}(p)$.
\item
First, if $p \in Hom_{+}(n,n)$ is {\bf H-Stable} then clearly the scaled polynomial $p_{{\bf c}}$ is also {\bf H-Stable}.\\
It follows from the definition of the {\bf Capacity} that $Cap(p_{{\bf c}}) = c_1 \cdots c_n Cap(p)$.\\
We get, by a direct computation, that 
$$
\frac{\partial^n}{\partial x_1\dots \partial x_n} p_{{\bf c}}(0,\dots,0) = c_1 \cdots c_n \frac{\partial^n}{\partial x_1\dots \partial x_n} p(0,\dots,0).
$$
This proves that the set of {\bf H-Stable} {\bf extremal} polynomials is invariant with respect to the scaling.
\end{enumerate}
\eprf

We need the following simple result (it was essentially proved in Lemma 3.8 from \cite{GS1}).
\fac \label{min}
Consider $p \in Hom_{++}(n,n), p(x_1,\dots x_n) = \sum_{r_1,\dots,r_n} a_{r_1,\dots,r_n} \prod_{1 \leq i \leq n} x_i^{r_{i}}$. Then there exists a positive vector 
${\bf t} =:(t_1,\dots,t_n) \in R^{n}_{++}, t_1 \cdots t_n =1$
such that
\beqn \label{exist}
p(t_1,\dots,t_n) = Cap(p) = \inf_{x_i > 0, 1 \leq i \leq n; x_1 \cdots x_n = 1} p(x_1,\dots,x_n).
\eeqn

Consider the corresponding scaled polynomial $p_{{\bf t}}$. Then the polynomial $q = \frac{p_{{\bf t}}}{p_{{\bf t}}(1,\dots,1)}$ is
doubly-stochastic. I.e. $\frac{\partial}{\partial x_i} q(1,1,\dots,1) = 1, 1 \leq i \leq n.$.
\efac

\prf
Consider a subset 
$$ 
T = \{(x_1,\dots,x_n) \in R^n_{++}: x_1 \cdots x_n = 1 \quad \mbox{and} \quad p(x_1,\dots,x_n) \leq p(1,\dots,1) \}.
$$
In order to prove that the infimum is attained, it is sufficient to show that the subset $T$ is compact.
Clearly, $T$ is closed, and we need to prove that $T$ is bounded. Let $(x_1,\dots,x_n) \in T$ and
assume WLOG that $\max_{1 \leq i \leq n} x_i = x_1$. Then
$$
a_{n,0,\dots,0} x_1^{n} \leq p(x_1,\dots,x_n) \leq p(1,\dots,1) \Rightarrow \max_{1 \leq i \leq n} x_i \leq \frac{p(1,\dots,1)}{a_{n,0,\dots,0}} < \infty.
$$
This shows the desired boundness of $T$ and the existence of the minimum.\\
Consider a positive vector $(t_1,\dots,t_n)$ such that $p(t_1,\dots,t_n) = Cap(p)$. Define $\alpha_i = \log(t_i), 1 \leq i \leq n$.
Then
$$
p(exp(\alpha_1),\dots,exp(\alpha_n)) = \min_{\beta_1+\dots \beta_n = 0} p(exp(\beta_1),\dots,exp(\beta_n)).
$$
Therefore there exists the Lagrange multiplier $\gamma$ such that
$$
\frac{\partial}{\partial \alpha_i}p(exp(\alpha_1),\dots,exp(\alpha_n)) = t_i \frac{\partial}{\partial t_i} p(t_1,\dots,t_n) = \gamma, 1 \leq i \leq n.
$$
It follows from the Euler's identity that $\gamma = p(t_1,\dots,t_n)$ and 
$$
\frac{\partial}{\partial x_i} q(1,\dots,1) = (p(t_1,\dots,t_n))^{-1} t_i \frac{\partial}{\partial t_i} p(t_1,\dots,t_n) = 1.
$$

\eprf
\rem
It is easy to prove that, in fact, the minimum in (\ref{exist}) is attained uniquely.
It was proved in \cite{stoc-06} that if $p \in Hom_{+}(n,n)$ is {\bf H-Stable} then
the minimum in (\ref{exist}) exists and attained uniquely iff
$$
\frac{\partial^n}{\partial x_j \partial x_j \prod_{m \neq (i,j)} \partial x_m} p(0,\dots,0),
\frac{\partial^n}{\partial x_i \partial x_i \prod_{m \neq (i,j)} \partial x_m} p(0,\dots,0) > 0: 1 \leq i \neq j \leq n.
$$
\erem

\prf ({\bf Proof of Theorem (\ref{uniq}).})\\
It follows from Proposition(\ref{ex-pr}) and Fact(\ref{min}) that we can assume
without loss of generality that the {\bf H-Stable} {\bf extremal} polynomial
$p \in Hom_{++}(n,n)$ is doubly-stochastic and all its coefficients
are positive real numbers.\\ 
\begin{enumerate}
\item{\bf Using uniqueness part of Lemma(\ref{derest})}\\
Let $\{e_1,\dots,e_n\}$ be the standard basis in $R^n$ and $p$ is now
a {\bf H-SStable} doubly-stochastic polynomial with positive
coefficients, $p$ satisfies the equality (\ref{EQ}). We need to look
at the case of equality in (\ref{Cap-In}). Recall the polynomial
$q_{n-1}$ is given by $q_{n-1}(x_1,\dots,x_{n-1}) = \frac{\partial}{\partial
x_n}p(x_1,\dots,x_{n-1},0)$.

Clearly, the polynomial $q_{n-1} \in Hom_{++}(n-1,n-1)$ also has positive coefficients. Let 
$$
q_{n-1}(t_{1,n},\dots,t_{n-1,n}) = \min_{x_1,\dots,x_{n-1} > 0, \prod_{1 \leq i \leq n-1} x_i =1} q_{n-1}(x_1,\dots,x_{n-1}).
$$
The existence of such a vector was proved in Proposition (\ref{min}).\\
It follows from the uniqueness
part of Lemma(\ref{derest}) that the univariate polynomial $R(t) =p(t_{1,n},\dots,t_{n-1,n},t) = p(\sum_{1 \leq i \leq n-1} t_{i,n-1} e_i + t e_n)$ has $n$ equal
negative roots: $R(t) = b(t+a_{n})^{n}; a_{n},b > 0$. This fact implies, as in the second part of Example (\ref{one-dim}),  that
$$
K_n =: e_n - \sum_{j \neq n} a_{j,n} e_j \in Null_{p}, a_{j,n} = \frac{t_{j,n}}{a_{n}} > 0.
$$
Since $p$ is doubly-stochastic,
hence it follows from (\ref{trace}) that the coordinates of $K_n$ sum to zero. Which gives that $\sum_{j \neq n} a_{j,n} = 1$. 

In the same way, we get that there exists an $n \times n$ column stochastic matrix $A$ with the zero diagonal and the positive off-diagonal part
such that the vectors
$$
K_i =: e_i - \sum_{j \neq i} a_{j,i} e_j \in Null_{p}, i \leq i \leq n.
$$
\item
Recall that our goal is to prove that $dim(Null_{p}) = n-1$. It follows that\\
$$
dim(Null_{p}) \geq dim(L(K_1,\dots,K_n)) = Rank(I-A),
$$
where $L(K_1,\dots,K_n)$ is the minimal linear subspace containing the set $\{K_1,\dots,K_n\}$.\\
Since the polynomial $p$ is non-zero thus $dim(Null_{p}) \leq n-1$. It is easy
to see that $Rank(I-A) =n-1$. Indeed, any principal $n-1 \times n-1$ submatrix of $I-A$
is {\it strictly diagonally dominant} and, therefore, is nonsingular.\\ 
We finally conclude that $dim(Null_{p}) = n-1$.

\end{enumerate}
\eprf

\section{Comments}
\begin{enumerate}
\item Falikman \cite{fal} and Egorychev \cite{ego} publications were followed by a flurry
of expository papers, which clarified and popularized the proofs. The author learned
the Egorychev's proof from \cite{knuth}. It is our guess that many scientists first learned about
Alexandrov inequalities for mixed discriminants and
Alexandrov-Fenchel inequalities for mixed volumes \cite{Al38} in one of those expository papers.
We would like to distinguish the following two papers: \cite{khov} and \cite{peetre}.
They both explicitly connected Alexandrov inequalities for mixed discriminants with
homogeneous hyperbolic polynomials. The paper \cite{khov} was, essentially a rediscovery
of Garding's theory \cite{gar}. Still, as the author had read \cite{khov} before reading \cite{gar},
the paper \cite{khov} gave us the first hint for the possibility of our approach.\\ 
The paper \cite{peetre}, apparently written as a technical report in 1981
and published only in 2006 in an obscure book, is technically very similar to \cite{khov}.
Besides, it implicitly introduced the Bapat's conjecture.\\
Other related publications are D. London's (univariate) papers \cite{london2},\cite{london3},\cite{london4}.\\

As far as we know, there were no previously published connections between
Shrijver-Valiant conjecture, which was thought to be of purely
combinatorial nature, and stable/hyperbolic polynomials.
\item
Two main ingredients of our approach, which make the proofs simple, are the
usage of the notion of {\bf Capacity} and Lemma (\ref{derest}). They together
allowed the simple induction. The induction, used in this paper, is by
partial differentiation. It is very similar to the inductive proofs of hyperbolic polynomials analogues of
Alexandrov inequalities for mixed discriminants
in \cite{khov},\cite{peetre}. Using our terminology, these analogues correspond
to the fact that the polynomial
$q_{2}(x_1,x_2) = \frac{\partial^{n-2}}{\partial x_{3}\dots \partial x_n}p(x_1,x_2,0,\dots,0)$ is either
zero or {\bf H-Stable} provided the polynomial $p$ is {\bf H-Stable}.\\
The idea to use {\bf Capacity} in the context of permanents is implicit in \cite{london1}.
The notion of {\bf Capacity} was crucial for algorithmic results in \cite{GS}, \cite{GS1} as
$\log \left(Cap(p) \right) = \inf_{y_1+\cdots +y_n =0} \log \left(p(e^{y_{1}},\dots,e^{y_{n}}) \right)$ and
the functional $\log \left(p(e^{y_{1}},\dots,e^{y_{n}}) \right)$ is convex for any polynomial with non-negative coefficients.\\
Probably, the papers \cite{GS}, \cite{GS1} were the first to reformulate Van der Waerden/Bapat conjectures
as in inequality (\ref{VDW}). Although quite simple, it happened to be a very enlighting observation.
\item
Our, inductive by the partial differentiation, approach was initiated in \cite{newhyp}. The main tool there was Vinnikov-Dubrovin determinantal representation \cite{vin}
of hyperbolic homogeneous polynomials in 3 variables. The paper \cite{newhyp} proved the implication\\
$Cap(p) > 0 \Longrightarrow  \frac{\partial^n}{\partial x_1\dots \partial x_n} p(0,\dots,0) > 0$ for
{\bf H-Stable} polynomials $p \in Hom_{+}(n,n)$. Additionally, it was proved that in this {\bf H-Stable} case  
the functional $Rank_{p}(S) = \max_{a_{r_1,\dots,r_n} \neq 0} \sum_{j \in S} r_j$ is submodular
and  
\beqn \label{rado}
a_{r_1,\dots,r_n} > 0 \Longleftrightarrow \sum_{j \in S} r_j \leq Rank_{p}(S): S \subset \{1,\dots,n\}.
\eeqn
The characterization (\ref{rado}) is a far reaching generalization of the Hall-Rado theorem.\\
The paper \cite{stoc-06} provides algorithmic applications of these results:
strongly polynomial deterministic algorithms for the membership problem as for the support
as well for the Newton polytope of {\bf H-Stable} polynomials $p \in Hom_{+}(m,n)$, given as oracles.
\end{enumerate}

\section{Acknowledgements}
The author is indebted to the anonymous reviewer for a very careful and thoughtful
reading of the original version of this paper. Her/his numerous corrections
and suggestions are reflected in the current version.\\
I would like to thank the U.S. DOE for financial support through
Los Alamos National Laboratory's {\bf LDRD} program.

\end{document}